

\magnification=\magstep 1
\def\mi{{\rm Mi}}

\def\i{{\cal I}}

\centerline{\bf The Eisenbud-Koh-Stillman Conjecture on Linear Syzygies}

\medskip\centerline{Mark L. Green\footnote*{Research partially supported by
the N.S.F.}, U.C.L.A.}\bigskip

Although the relationship between minimal free resolutions and Koszul
cohomology
has been known for a long time, it has been difficult to find a way to fully
utilize the ``exterior" nature of the Koszul classes.   The technique used here
seems to be one way to begin to do this.  We prove a conjecture of
Eisenbud-Koh-Stillman on linear syzygies and in consequence a conjecture of
Lazarsfeld and myself on points in projective space.  The main novelties in the
proof are the use of ``exterior minors," explained below, and showing that
certain kinds of linear syzygies in the exterior algebra are impossible.
I will work over a field of arbitrary characteristic.

It is a pleasure to acknowledge David Eisenbud for many highly useful
conversations regarding this work.  In particular, he simplified several
arguments in the original version, some of which had only worked in
characteristic 0.

\proclaim DEFINITION.  Consider two vector spaces $A$,$B$ of dimensions $a$,$b$
respectively, and let $V$ be a vector space of dimension $n$.  Consider a
$b\times a$ matrix of linear forms, which we think of as a linear map $M\colon
A\to B\otimes V$.  By a {\bf generalized column} of $M$ we mean, for some
non-zero $\alpha\in A$, the map $M(\alpha)\colon B^*\to V$, and by the {\bf
rank
of a generalized column} $\alpha$ we mean the rank of the map $M(\alpha)$;
similarly an element $\beta^*\in B^*$ gives a {\bf generalized row} which is a
map $M(\beta^*)\colon A\to V$ whose rank is the rank of $M(\beta^*)$.  Now
$\wedge^kM\colon \wedge^kA\to \wedge^k(B\otimes V)$.  There are natural maps
$P_s\colon \wedge^k(B\otimes V)\to \wedge^kB\otimes S^kV$ and $P_e\colon
\wedge^k(B\otimes V)\to S^kB\otimes\wedge^kV$.  The maps
$P_s\circ\wedge^kM\colon \wedge^kA\to \wedge^kB\otimes S^kV$ and
$P_e\circ\wedge^kM\colon\wedge^kA\to S^kB\otimes \wedge^kV$ induce natural maps
$\mi^k_s(M)\colon\wedge^kA\otimes\wedge^kB^*\to S^kV$ and $\mi^k_e(M)\colon
\wedge^kA\otimes S^kB^*\to \wedge^kV$.  The image of $\mi^k_s(M)$ is just the
usual ideal $I^k_s(M)$ of {\bf k by k minors} of $M$; the image of $\mi^k_e(M)$
we will denote by $I^k_e(M)$ and will call the {\bf k by k exterior minors} of
$M$.

The exterior minors are quite interesting and useful.  I do not know of a good
reference for their properties, so I will prove what I need.

\proclaim PROPOSITION 1.  Let $M$ be a $b\times a$ matrix of linear forms such
that every generalized column of $M$ has rank $b$.  Then the map $\mi^a_e(M)$
is
injective, and hence there are $b-1+a\choose a$ linearly independent
$a\times a$ exterior minors of $M$.

\noindent PROOF:  For a non-zero $\alpha\in A$, the map $M(\alpha)\colon B^*\to
V$ has rank $b$, and thus the set of subspaces $W\subseteq V$ of codimension
$a+b-1$ which meet the image of $M(\alpha)$ has codimension $\ge a$ in the
Grassmannian.  It follows that the set of subspaces $W$ of this dimension
meeting the image of some $M(\alpha)$ as $\alpha$ ranges over ${\bf P}(A)$
has codimension $\ge 1$.  If we replace $M$ by the composition $M'\colon A\to
B\otimes V\to B\otimes V'$ obtained from a general projection $V\to V'$ to a
vector space of dimension $a+b-1$, the hypothesis continues to hold, and the
exterior minors of $M'$ are the projection of the exterior minors of $M$ under
the map $\wedge^a V\to \wedge^a V'$.  It is thus enough to treat the case ${\rm
dim}(V)= a+b-1$.

We may regard $M$ as an $(a+b-1)\times a$ matrix of linear forms in $B$.  The
exterior minors of $M$ are the usual minors of this new matrix.  The hypothesis
on generalized columns of $M$ translates into the hypothesis that the new
matrix never drops rank.  The Eagon-Northcott complex (cf. Appendix A2 of [6])
now shows that these minors are linearly independent.  This completes the
proof.

The way we will make use of exterior minors is the following result, which is
essentially the same as what happens in the commutative case:

\proclaim PROPOSITION 2.  If $M\colon A\to
B\otimes V$ is a linear map and $f_M\colon \wedge^m V\otimes A \to
\wedge^{m+1}V\otimes B$ is the naturally associated map, then $I^a_e(M)$
annihilates ${\rm ker}(f_M)$, i.e. if $\phi\in I^a_e(M)$ and $\alpha\in {\rm
ker}(f_M)$, then $\phi\wedge\alpha\in \wedge^{m+a}V\otimes A$ is zero.

\noindent PROOF:  Let $a_1,\ldots ,a_a$ be a basis for $A$, and let $\alpha
=\sum_{i=1}^a \alpha_i\otimes a_i$ with $\alpha_i\in \wedge^mV$ for all $i$.
Let $M(a_i)\in B\otimes V$ be the image of $a_i$ under $M$.  The hypothesis is
that $\alpha\in {\rm ker}(f_M)$ is equivalent to $\sum_i \alpha_i\wedge M(a_i)
=0$, where by $\wedge$ we mean that we multiply elements of $B$ symmetrically
and
elements of $V$ anti-symmetrically, with the result that $\wedge$
anti-commutes.  What we need to show is that for all $i$,
$M(a_1)\wedge\cdots\wedge M(a_n)\wedge\alpha_i =0$.  However, we may write the
left-hand side as $M(a_1)\wedge\cdots\wedge M(a_{i-1})\wedge(\sum_j
M(a_j)\wedge\alpha_j)\wedge M(a_{i+1})\wedge\cdots\wedge M(a_n)$, and this is
zero.

\smallskip\noindent {\bf REMARK.} A more elegant approach, suggested by
Eisenbud,
is to notice that ${\rm Mi}_e^{a-1}(M)$ gives a map $B\otimes V\to S^aB\otimes
\wedge^{a-1}A^*\otimes\wedge^a V$, which then maps naturally to
$S^aB\otimes \wedge^aA^*\otimes \wedge^aV\otimes A$.  The composition is a map
$M'\colon B\otimes V\to S^a B\otimes \wedge^aA^*\otimes\wedge^aV\otimes A$
which functions as a ``companion matrix" to $M$ because the composition $M'M=
{\rm Mi}_e^a(M)\otimes {\rm id}_A$.  This formula implies Proposition 2.

\proclaim PROPOSITION 3. Let $V$ be a vector space of dimension $n$ and
$W\subseteq \wedge^{n-p}V$ a linear subspace of dimension $p>0$. Then there
exists a $0<k\le p$ and a $(p-k)$-dimensional subvariety $Z\subseteq G(p-1,W)$
of
$(p-1)$-dimensional subspaces $U\subset W$ such that for all $U\in Z$, the
image
of $U\otimes V\to \wedge^{n-p+1}V$ has codimension $\ge k$ in the image of
$W\otimes V\to \wedge^{n-p+1}V$.

\noindent PROOF:  Let $B^*={\rm ker}(W\otimes V\to \wedge^{n-p+1}V)$.  Let $M$
be
the $b\times p$ matrix of linear relations in $\wedge^*V$ of $W$, which we view
as a map $W^*\to B\otimes V$.  Let $I_{W,V}$ and $I_{U,V}$ be the images of
$W\otimes V$ and $U\otimes V$ respectively in $\wedge^{n-p+1}V$.  Let $C^* =
{\rm ker}(U\otimes V\to \wedge^{n-p+1}V)$.  There is then an exact sequence
$$0\to {B^*\over C^*}\to (W/U)\otimes V\to {I_{W,V}\over I_{U,V}}\to 0.$$We
conclude that ${\rm dim}(I_{W,V}/I_{U,V}) = n-{\rm dim}(B^*/C^*)$.  If $w^*\in
W^*$ is the annihilator of $U$, then ${\rm dim}(B^*/C^*)$ is just the rank of
the generalized column of $M$ corresponding to $w^*$.  In the ${\bf
P}^{p-1}$ parametrizing generalized columns of $M$, let $Z_r$ be the space of
generalized columns of rank $r$, and $d(r)={\rm dim}(Z_r)$.  The negation of
the
conclusion of the proposition is that $Z_r =\phi$ for $r\le n-p$ and
$d(r)<r+p-n$ for all $r>n-p$.  We will assume this and derive a contradiction.

Let $\bar M$ be a $(n-p+1)\times p$ matrix obtained by
choosing $(n-p+1)$ general generalized rows of $M$.  For any given
generalized column of $M$, let $r$ be its rank.  If $r\ge n-p+1$, which we
may assume, then the set of projections of $M$ to a $p\times (n-p+1)$ matrix
for
which this generalized column does not have maximal rank has codimension
$r+p-n$.  Every generalized column of $\bar M$
has maximal rank provided $d(r)<r+p-n$ for all $r$, as then a general
projection does not belong to the ``bad set" for any generalized column.

Thus a general choice of $\bar M$ has every generalized column of rank
$n-p+1$.  By the first proposition, there are at least $n\choose p$ linearly
independent elements of $I^p_e(\bar M)\subseteq \wedge^pV$, and hence
$I^p_e(\bar M) =\wedge^pV$.  By the second proposition, the elements of
$I^p_e(\bar M)$ annihilate $W$, and therefore $W=0$, which is a contradiction.

As a consequence of the foregoing proposition, we obtain the following result
and its corollary, which were conjectured by Eisenbud, Koh, and Stillman (see
[1].) We use the notation ${\cal K}_{p,q}(M,V)$ to denote the Koszul cohomology
group $H^p(\wedge^\bullet V\otimes M_{q+p-\bullet})$.

\proclaim THEOREM 4.  Let $M=\oplus_{q\ge 0}M_q$ be a finitely generated
$S(V)$-module and assume ${\rm rank}(M_0) = p > 0$.  Let $R\subseteq M_0\otimes
V$ be the module of relations.  Then if the (affine) dimension of the rank one
relations $R_1$ has dimension $<p$, then the Koszul cohomology group $${\cal
K}_{k,0}(M,V) ={\rm ker}(\wedge^kV\otimes M_0\to \wedge^{k-1}V\otimes
M_1)$$vanishes for all $k\ge p$.

\noindent PROOF:  We may proceed by induction on $p$, the case $p=1$ being
obvious. As is well-known, since $M_q=0$ for $q<0$, the vanishing of ${\cal
K}_{p,0}(M,V)$ would imply the vanishing of ${\cal K}_{k,0}(M,V)$ for all $k\ge
p$.  Let $\alpha\in\wedge^pV\otimes M_0$ be a non-zero element of ${\cal
K}_{p,0}(M,V)$.  Under the map $\wedge^pV\otimes M_0\to \wedge^{p-1}V\otimes
V\otimes M_0$, $\alpha$ must map to an element of $\wedge^{p-1}V\otimes
R$.  Let $W$ be the image of the map $M_0^*\to \wedge^pV$ given by $\alpha$. We
may assume that it has dimension $p$, since otherwise we could shrink $M_0$ and
$p$. For any $\beta\in\wedge^{p-1}V^*$, the contraction
$<\beta,\alpha>\in V\otimes M_0$ automatically lies in $R$. Any $m\in M_0$
annihilates a $(p-1)$-dimensional linear subspace of $M_0^*$; let $U_m$ be its
image under the map $M_0^*\to \wedge^pV$ determined by $\alpha$. If the image
of
$U_m\otimes V^*\to\wedge^{p-1}V$ has codimension $k_m$ in the image of
$W\otimes
V^*\to \wedge^{p-1}V$, then we obtain a $k_m$-dimensional linear space of
non-trivial rank one relations in $M_0\otimes V$ lying in $m\otimes V$.

Now choose a generator $\tau$ for $\wedge^nV^*$, and let $W' =
<\tau,W>\subseteq\wedge^{n-p}V^*$ and $U_m' =<\tau ,U_m>$.  Then $k_m$ is
also the codimension of the image of $U_m'\otimes V^*\to \wedge^{n-p+1}V^*$
in the image of $W'\otimes V^*\to \wedge^{n-p+1}V^*$.  We now invoke the
preceding proposition to conclude that for some $0<k\le p$, there is a
variety in $G(p-1,M_0)$ of dimension at least $p-k$ such that for all
$m$ in this variety, there exists a $k$-dimensional family of rank one
relations
of $M$ in $m\otimes V$.  This completes the proof.

\proclaim DEFINITION. A relation of rank $\le r$ is a non-zero
element of ${\rm ker}(S\otimes V\to M_1)$ for some linear subspace $S\subseteq
M_0$ of rank $r$.  We will say that such a relation {\bf involves the
linear subspace} $S'\subset M_0$ if $S'\subseteq S$.  A corollary of the
Theorem is:

\proclaim COROLLARY 5.  Let $M=\oplus_{q\ge 0}M_q$ be a finitely generated
$S(V)$-module and let ${\rm rank}(M_0)= m_0>0$. If ${\cal K}_{p,0}(M,V)\ne 0$,
then for a general choice of $(m_0-p)$-dimensional subspace $S\subseteq M_0$,
the affine dimension of the rank $\le (m_0-p+1)$ relations involving $S$ is at
least $p$.

\noindent PROOF: If $S\subseteq M_0$ is a linear subspace, $\tilde S$ the
submodule of $M$ it generates, and let $\bar M = M/ \tilde S$.  For a general
choice of $S$ of dimension $m_0-p$, our Koszul class $\alpha\in \wedge^p
V\otimes M_0$ maps to a non-zero class in $\wedge^p V\otimes (M_0/S)$.  Rank 1
relations for $\bar M$ lift to rank $\le m_0-p+1$ relations involving $S$.
However, ${\rm rank}(\bar M_0) = p$, so the preceding theorem applies,
completing the proof.

\smallskip\noindent {\bf REMARK.} Eisenbud, Koh, and Stillman noted when they
made their conjecture that Corollary 5 would follow from Theorem 4.  I suspect
that one can improve this further using mixed minors (see below.)

After Eisenbud, Koh, and Stillman made their conjecture, I showed [2]
that it implies part of the following conjecture [4] of Rob Lazarsfeld and
myself.  However, in fact one can get the full conjecture:

\proclaim THEOREM 6.  Let $Z$ be a set of $2r+1-p$ points in ${\bf P}^r$, $p\ge
0$.  Then either $Z$ has property $N_p$ or there exists a subset $Z'\subseteq
Z$
and a linear space $L$ such that $Z'\subseteq L$, $Z'$ consists of at least $2
{\rm \ dim}(L)+2-p$ points and property $N_p$ fails for $Z'$.

For the definition of property $N_p$, see [5].

\smallskip\noindent PROOF:Let $\i_Z$ denote the ideal
sheaf of $Z$ and $V= H^0({\cal O}_{{\bf P}^r}(1))$.  From the
surjective map $H^0({\cal O}_Z(k))\to H^1(\i_Z(k))$, we conclude that there is
an injection $0\to H^1(\i_Z(1))^*\to H^0({\cal O}_Z(1))^*$.  If we let
$Z=\{P_1,\ldots ,P_{2r+1-p}\}$, and let $v_i$ denote the element of $H^0({\cal
O}_Z(k))^*$ given by evaluation at $P_i$ (this depends on a choice of
trivialization of the hyperplane bundle at each $P_i$), then the map $V\otimes
H^1(\i_Z(k))^*\to H^1(\i_Z(k-1))^*$ is given by $$l\otimes \sum_i a_iv_i\mapsto
\sum_i l(P_i)a_iv_i.$$

The case $p=0$ (done with Rob Lazarsfeld at the time of [4])
proceeds as follow---choose a subscheme $Z'$  of $Z$  which violates $N_0$,
but such that no proper subscheme of $Z'$ violates $N_0$.  This implies that if
we write a non-zero element $\phi\in H^1({\cal I}_{Z'}(2))^*$ as $\phi
=\sum_{i=1}^{2r+1} \phi_i v_i$, then $\phi_i\ne 0$ for all $i$.  Since
$h^1({\cal I}_{Z'}(1))\le  r <r+1$, under the multiplication $V\otimes
H^1({\cal
I}_{Z'}(2))^*\to H^1({\cal I}_{Z'}(1))^*$, there is some linear form $h\in V$
annihilating $\phi$.  Thus $\sum_i h(P_i) \phi_i v_i =0$.  Since the $v_i$ are
linearly independent, this implies that $h(P_i)=0$ for all $i$, and thus $Z'$
lies on the hyperplane $h$.  Now  either $Z'$ consists of $\le 2(r-1)+1$
points,
in which case we proceed inductively on the dimension, or it has $\ge 2(r-1)+2$
points, in which case we are done.  Of course, if $H^1({\cal I}_Z(2))=0$, then
${\cal I}_Z$ is 3-regular.  Thus if property $N_p$ fails for $p>0$, it must be
that ${\cal K}_{k,3}(I_Z,V)\ne 0$ for some $k\le p-1$, and we may reduce to the
case $k=p-1$.

   From the Koszul complex
$$0\to \wedge^{r+1}V\otimes \i_Z(-r-1)\to\cdots\to V\otimes \i_Z(-1)\to \i_Z\to
0$$ twisted by ${\cal O}_{{\bf P}^r}(p+3)$, we see that the Koszul group ${\cal
K}_{p-1,3}(I_Z,V)\cong {\cal K}_{p+1,1}(M,V)$, where $M=\oplus_kM_k
=\oplus_kH^1(\i_Z(k))$.  Since ${\cal K}_{p,2}(Z,V)\cong {\cal
K}_{p-1,3}(I_Z,V)$, we say that this group being non-zero implies that ${\cal
K}_{r-p,1}(M^*,V)\ne 0$, where $M^*_k = H^1(\i_Z (k))^*$ (in this module,
multiplication decreases degree.)  We may without loss of generality assume
that
$Z$ is not contained in a hyperplane.  Thus $h^1(\i_Z(1))= r-p$, from which we
conclude by the Theorem that for any non-zero class $\lambda\in {\cal
K}_{r-p,1}(M^*,V)$, there are at least an $(r-p)$-dimensional family of rank 1
relations in ${\rm im}(\wedge^{r-p-1}V^*\to V\otimes H^1(\i_Z(1))^*)$, where
the
map is induced by $\lambda$.  The $v_i$ are linearly independent, and thus a
relation of rank 1 is an element $l\otimes \sum_i a_iv_i$ such that, for all
$i$, either $l(P_i)=0$ or $a_i=0$.  On any irreducible component of the
intersection of the image of $\lambda$ in $V\otimes H^1(\i_Z(1))^*$ with the
rank one locus, there is a decomposition $Z= Z' +Z''$ of $Z$ into disjoint
subsets such that $l(P_i)=0$ for $P_i\in Z'$ and $a_i=0$ if $P_i\in Z''$.
{}From
the exact sequence $$0\to H^1(\i_{Z'}(1))^*\to H^1(\i_Z(1))^*\to H^0({\cal
O}_{Z''}(1))^*$$we see that the $\sum_i a_iv_i$ actually belong to
$H^1(\i_{Z'}(1))^*$ in this circumstance, so that the image of $\lambda$
intersect the rank one locus gives a subvariety $Y\subseteq
H^0(\i_{Z'}(1))\times H^1(\i_{Z'}(1))^*$ of dimension $\ge r-p$ consisting of
rank 1 relations for $M^*$.  This implies that
$h^0(\i_{Z'}(1))+h^1(\i_{Z'}(1))\ge r-p$.  If $Z'$ is $m$ points spanning a
${\bf P}^k$, then  $h^0(\i_{Z'}(1))=r-k$ and $h^1(\i_{Z'}(1))= m-1-k$, so that
the preceding inequality becomes $m\ge 2 k+ 1-p$.   Possibly by enlarging $Z'$,
we may without loss of generality in this construction assume that for all
$p\in
Z''$, there exists an $(l,\phi)\in Y$ such that $l(p)\ne 0$.

It remains to show that property $N_p$ fails for $Z'$.  Let
$\phi_1,\ldots, \phi_m, \psi_1,\ldots, \psi_t$ be a basis for $H^1(\i_Z(1))^*$
such that the $\phi_i$ are a basis for $H^1(\i_{Z'}(1))^*$.  If we write our
Koszul class as $\lambda =\sum_i \alpha_i\otimes \phi_i +\sum_j
\beta_j\otimes \psi_j$, then if $x_1,\ldots x_{r+1}$ is a basis for $V$  and
$e_1,\ldots ,e_{r+1}$ is the dual basis, then the condition that $\lambda$ be a
Koszul class is $$\sum_{i,\nu} <\alpha_i,e_\nu>\otimes x_\nu\phi_i
+\sum_{j,\nu}<\beta_j,e_\nu>\otimes x_\nu\psi_j =0.$$We read off that the image
of $\sum_{j,\nu}<\beta_j,e_\nu>\otimes x_\nu\psi_j =0$ in
$\wedge^{p-1}V\otimes H^0({\cal O}_{Z''}(2))^*$.  However, the $\psi_j$ are
linearly independent there, and thus $\sum_\nu <\beta_j,e_\nu>\otimes x_\nu=0$
on
$Z''$ for all $j$. It follows that $<\beta_j,v>=0$ for all $v\in \bar Z''$ and
all $j$, where $\bar Z''$ is the linear span of $Z''$.  It thus follows that
$<\lambda, v> =\sum_i <\alpha_i,v>\otimes \phi_i$ if $v\in \bar Z''$, and
$<\lambda,v>$ thus determines an element of ${\cal K}_{r-p-1,1}(\bar
M^*,v^\perp)$, where $\bar M =\oplus_{q\ge 0} H^1({\cal I}_{Z'}(q))^*$.  Thus
$Z'$ fails to have property $N_p$ when viewed as a subset of $v^\perp$ unless
$<\lambda,v>=0$ in ${\cal K}_{r-p-1,1}(\bar M,v^\perp)$ for all $v\in \bar
Z''$.  Note that $H^1(\i_{Z'}(2))^*=0$ unless $Z'$ fails to have property $N_0$
and hence {\sl a fortiori} $N_p$.  We have arranged things so that the image of
$\lambda$ intersected with the rank 1 locus does not lie in $(\bar
Z'')^\perp\otimes H^1(\i_{Z'}(1))^*$.  The proof is now done, modulo the
following elementary lemma:

\proclaim LEMMA 7. Let $M$ be a finitely generated module over $S(V)$.  If
$\lambda\in {\cal K}_{p,q}(M,V)$ and $v\in V^*$, then the contraction
$<\lambda, v>\in {\cal K}_{p-1,q}(M,v^\perp )$.  If $M_{q-1}=0$ and the image
of
the map $\wedge^{p-1}V^*\to V\otimes M_q$ induced by $\lambda$ is not contained
in $v^\perp\otimes M_q$, then $<\lambda,v>\ne 0$ in ${\cal K}_{p-1,q}(M,
v^\perp )$.

\noindent PROOF:  If $m_1,\ldots ,m_k$ is a basis for $M_q$, then if
$\lambda=\sum_i\lambda_i\otimes m_i$ with $\lambda_i\in \wedge^pV$, the Koszul
condition is $\sum_{i,\nu} <\lambda_i,e_\nu>\otimes x_\nu m_i=0$; here
$x_1,\ldots ,x_n$ is a basis for $V$ and $e_1,\ldots ,e_n$ the dual basis for
$V^*$.  We may take $v=e_1$ and $v^\perp = {\rm span}(x_2,\ldots ,x_n)$.  If we
contract $e_1$ with the equality above, we get $\sum_{i=1}^k\sum_{\nu=1}^n
<\lambda_i, e_1\wedge e_\nu>\otimes x_\nu m_i =0$.  Of course, we may sum over
$2\le \nu\le n$ in this equality, and then we have the Koszul condition for
$<\lambda,e_1>$ over $v^\perp$.  Further, the image of the map induced by
$\lambda$ does not belong to $v^\perp\otimes M_q$ if and only if $\sum_i
<\lambda_i,e_1>\otimes m_i \ne 0$, and this is equivalent to
$<\lambda,e_1>\ne 0$.  If $M_{q-1}=0$, then this implies that $<\lambda,e_1>$
does not vanish in  ${\cal K}_{p-1,q}(M, v^\perp )$.

\smallskip\noindent {\bf REMARK.} If $Z$ is $2r+2-p$ points on a rational
normal curve in ${\bf P}^r$, then property $N_p$ fails for $Z$.

\proclaim REMARK.  A linear map $M\colon A\to B\otimes V$ gives maps
$\wedge^kA\to B^\lambda\otimes V^{\lambda^t}$ for $\lambda$ any Young diagram
of
size $k$, using the Cauchy decomposition $$\wedge^k(B\otimes V)\cong
\oplus_\lambda  (B^\lambda\otimes V^{\lambda^t}).$$One can do similar things
using the decomposition of $(B\otimes V)^\lambda$ for other Young diagrams
$\lambda$.  I call these the {\bf mixed minors} of the matrix $M$ of linear
forms.  I hope to give some applications of mixed minors in a later paper.

\medskip\noindent\underbar{BIBLIOGRAPHY}\medskip\frenchspacing

\item{[1]} D. Eisenbud and J. Koh, ``Some linear syzygy conjectures," Adv. in
Math. 90 (1991), 47-76.

\item{[2]} D. Eisenbud and J. Koh, ``Remarks on points in projective space,"
in {\sl Commutative Algebra}, Conf. Proc. of the 1987 conference in Berkeley
(M. Hochster, C. Huneke, and J.D. Sally, Eds.), 157-172, MSRI Publ 15,
Springer-Verlag, New York, 1989.

\item{[3]} M. Green and R. Lazarsfeld, ``On the projective normality
of complete linear
series on an algebraic curve," Inv. Math. 83 (1986), 73-90.

\item{[4]} M. Green and R. Lazarsfeld, ``Some results on the syzygies of
finite sets and algebraic curves," Comp. Math. 67 (1988), 301-314.

\item{[5]} R. Lazarsfeld, ``A sampling of vector bundle techniques in the study
of linear series," in {\sl Lectures on Riemann Surfaces}, World Scientific,
Singapore (1989), 500-559.

\item{[6]} D. Eisenbud, {\sl Commutative Algebra with a View Toward Algebraic
Geometry}, Spring\-er-Verlag, New York (1995).

 \end